\documentclass{amsart}
\usepackage{longtable}
\usepackage{graphicx}
\usepackage{amsmath}
\usepackage{amssymb}
\newtheorem{theorem}{Theorem}[section]
\newtheorem{proposition}[theorem]{Proposition}

\newcommand{\tabtopsp}[1]{\vbox{\vbox to#1{}\vbox to1ex{}}}

\theoremstyle{definition}

\newcommand{\F}{\Bbb F}

\theoremstyle{remark}
\newtheorem{remark}[theorem]{Remark}

\begin{document}

\title{A partial order on the set of 
prime knots with up to 11 crossings}

\author{
Keiichi Horie, 
Teruaki Kitano, 
Mineko Matsumoto and 
Masaaki Suzuki} 

\date{\today}

\address{Software Cradle Co., Ltd., 
3-4-5, Umeda, Kita-ku Osaka 530-0001, Japan}

\email{horie@cradle.co.jp}

\address{Department of Information Systems Science,
Faculty of Engineering,
Soka University, 
1-236 Tangi-cho, Hachioji-shi, 
Tokyo, 192-8577, Japan}

\email{kitano@soka.ac.jp}

\address{Department of Information Systems Science,
Faculty of Engineering,
Soka University
1-236 Tangi-cho, Hachioji-shi, 
Tokyo, 192-8577, Japan}

\email{e08m5223@soka.ac.jp}

\address{Department of Mathematics,
Akita University, 
1-1 TegataGakuenmachi, 
Akita, 010-8502, Japan}

\email{macky@math.akita-u.ac.jp}



\thanks{The second author is supported
in part by Grand-in-Aid for Scientific Research
(No. 17540064),
The Ministry of Education, Culture,
Sports, Science and Technology, Japan.}

\begin{abstract}
Let $K$ be a prime knot in $S^3$ and $G(K)=\pi_1(S^3-K)$ the knot group. 
We write $K_1 \geq K_2$ if there exists a surjective homomorphism 
from $G(K_1)$ onto $G(K_2)$. 
In this paper, we determine 
this partial order on the set of prime knots with up to $11$ crossings. 
There exist such $801$ prime knots and then $640,800$ should be considered. 
The existence of a surjective homomorphism can be proved by constructing it explicitly. 
On the other hand, the non-existence of a surjective homomorphism can be proved 
by the Alexander polynomial and the twisted Alexander polynomial. 
This work is an extension of the result of \cite{KS1}. 
\end{abstract}

\maketitle

\section{Introduction}\label{intro}

Let $K$ be a prime knot in $S^3$ and $G(K)=\pi_1(S^3-K)$ the knot group. 
For two prime knots $K_1,K_2$, we write $K_1 \geq K_2$ 
if there exists a surjective homomorphism from $G(K_1)$ onto $G(K_2)$. 
It is well known that 
this relation is a partial order on the set of prime knots. 
See \cite{Go, Gi} as references. 
 
Some geometric reasons for the existence of a surjective homomorphism are known, 
that is, degree one map and periodic knot. 
See \cite{Go, Gi, Wang1} for example. 
It could be considered that this partial order measures 
some kind of complexity for knots. 
Then it should be determined 
whether there exists a surjective homomorphism between the knot groups of two given knots. 

The first useful criterion is given by the Alexander polynomial. 
A well-known property on the Alexander polynomial gives us a necessary condition 
for the existence of a surjective homomorphism. 
However, it is not sufficient to decide the non-existence of a surjective homomorphism 
between knot groups. 
In 1990's, 
the Alexander polynomial was generalized to the twisted Alexander polynomial 
for a knot with a linear representation in \cite{L, Wa1, JW}. 
The above necessary condition 
can also be extended to that of the twisted Alexander polynomial 
in \cite{KSW}, 
which is more effective 
for determining the non-existence of a surjective homomorphism. 

By using the above criteria, 
two of the authors \cite{KS1} completely determined this partial order 
for the Rolfsen's knot table \cite{R}, 
which contains $249$ prime knots with up to $10$ crossings. 
There are $61,752$ cases to be treated. 

\begin{theorem}[\cite{KS1}]
The partial order ``$\geq$'' on the set of prime knots with up to $10$ crossings 
is given by 
\[
\left.
\begin{array}{l}
8_5,8_{10},8_{15},8_{18},8_{19},8_{20},8_{21},\\
9_1,9_6,9_{16},9_{23},9_{24},9_{28},9_{40},\\
10_5,10_9,10_{32},10_{40},10_{61},10_{62},10_{63},10_{64},10_{65},10_{66},\\
10_{76},10_{77},10_{78},10_{82},10_{84},10_{85},10_{87},10_{98},10_{99},10_{103},\\
10_{106},10_{112},10_{114},10_{139},10_{140},10_{141},10_{142},10_{143},10_{144},10_{159},10_{164}
\end{array}
\right\}\geq 3_1 ,
\]
\[
8_{18},9_{37},9_{40},
10_{58},10_{59},10_{60},10_{122},10_{136},10_{137},10_{138}
\geq 4_1 ,
\]
\[
10_{74},10_{120},10_{122} \geq 5_2 .
\]
\end{theorem}

The purpose of this paper is to  determine 
this partial order ``$\geq$'' on the set of prime knots with up to $11$ crossings, 
because some phenomena require more higher crossing knots. 
There exist $801$ prime knots with up to $11$ crossings 
and then $640,800$ cases should be investigated. 

The following is the main result of this paper. 

\begin{theorem}\label{mainthm} 
The partial order ``$\geq$'' on the set of prime knots with up to $11$ crossings 
is given by 
\[
\left.
\begin{array}{l}
8_5,8_{10},8_{15},8_{18},8_{19},8_{20},8_{21},\\
9_1,9_6,9_{16},9_{23},9_{24},9_{28},9_{40},\\
10_5,10_9,10_{32},10_{40},10_{61},10_{62},10_{63},10_{64},10_{65},10_{66},\\
10_{76},10_{77},10_{78},10_{82},10_{84},10_{85},10_{87},10_{98},10_{99},10_{103},\\
10_{106},10_{112},10_{114},10_{139},10_{140},10_{141},10_{142},10_{143},10_{144},10_{159},10_{164}, \\
11a_{43}, 11a_{44}, 11a_{46}, 11a_{47}, 11a_{57}, 11a_{58}, 
11a_{71}, 11a_{72}, 11a_{73}, 11a_{100}, \\
11a_{106}, 11a_{107}, 11a_{108}, 11a_{109}, 11a_{117}, 
11a_{134}, 11a_{139}, 11a_{157}, 11a_{165}, 11a_{171}, \\
11a_{175}, 11a_{176}, 11a_{194}, 11a_{196}, 11a_{203}, 
11a_{212}, 11a_{216}, 11a_{223}, 11a_{231}, 11a_{232}, \\
11a_{236}, 11a_{244}, 11a_{245}, 11a_{261}, 11a_{263}, 
11a_{264}, 11a_{286}, 11a_{305}, 11a_{306}, \\
11a_{318}, 11a_{332}, 11a_{338}, 11a_{340}, 11a_{351}, 
11a_{352}, 11a_{355}, 11n_{71}, 11n_{72}, 11n_{73},\\
11n_{74}, 11n_{75}, 11n_{76}, 11n_{77}, 11n_{78}, 11n_{81}, 
11n_{85}, 11n_{86}, 11n_{87}, 11n_{94}, \\
11n_{104}, 11n_{105}, 11n_{106}, 11n_{107}, 11n_{136}, 
11n_{164}, 11n_{183}, 11n_{184}, 11n_{185} 
\end{array}
\right\}\geq 3_1 ,
\]
\[
\left.
\begin{array}{l}
8_{18},9_{37},9_{40},\\
10_{58},10_{59},10_{60},10_{122},10_{136},10_{137},10_{138},\\
11a_{5}, 11a_{6}, 11a_{51}, 11a_{132}, 11a_{239}, 11a_{297}, 11a_{348}, 11a_{349}, \\
11n_{100}, 11n_{148}, 11n_{157}, 11n_{165}
\end{array}
\right\}
\geq 4_1 ,
\]
\[
11n_{78}, 11n_{148} \geq 5_1 ,
\]
\[
10_{74},10_{120},10_{122}, 11n_{71}, 11n_{185}
\geq 5_2 ,
\]
\[11a_{352} \geq 6_1,
\]
\[11a_{351} \geq 6_2 ,
\]
\[
11a_{47}, 11a_{239} \geq 6_3 .
\]
\end{theorem}

\begin{remark}
In this paper, 
the numbering of the knots 
with $11$ crossings follows 
that of the web page ``KnotInfo'' \cite{LC} 
by Cha and Livingston. 
\end{remark}

In Section \ref{recipe}, 
we explain the recipe to obtain the main result. 
In Section \ref{problems}, we put some problems. 
In Section \ref{table}, 
some tables used to prove Theorem \ref{mainthm} are described. 
Tables \ref{surj-31}-\ref{surj-63} are lists of surjective homomorphisms 
and Table \ref{non-ex} lists data to check the non-existence. 

\section{Proof}\label{recipe}

In this section, 
we shall explain how to obtain Theorem \ref{mainthm}. 
It is sufficient to prove Theorem \ref{mainthm} for any pair of knots 
where at least one of two knots  has  11 crossings, 
because we have already proved in \cite{KS1} 
for all pairs of knots with up to $10$ crossings. 
Therefore the number of the cases to be considered is $579,084$. 
We prove Proposition 2.1 and 2.2 for the main result. 

\begin{proposition}\label{prop-ex}
There exist surjective group homomorphisms 
$G(K_1) \to G(K_2)$ for all pairs of knots in Theorem \ref{mainthm}. 
\end{proposition}

We can construct explicitly surjective homomorphisms as shown 
in Tables \ref{surj-31}-\ref{surj-63}. 
Next we have to prove the following. 

\begin{proposition}\label{prop-nonex}
There does not exist a surjective group homomorphism 
$G(K_1) \to G(K_2)$ for any pair of knots that does not appear 
in Theorem \ref{mainthm}. 
\end{proposition}

To prove Proposition \ref{prop-nonex}, 
we consider the Alexander polynomial 
and the twisted Alexander polynomial of a knot. 

First, recall the following well-known proposition 
on the Alexander polynomial $\Delta_K(t)$. 
See \cite{CF} as a reference.  

\begin{proposition}\label{alex-div}
If there  exists a surjective homomorphism 
from $G(K_1)$ onto $G(K_2)$, 
then $\Delta_{K_1}(t)$ is divisible by $\Delta_{K_2}(t)$. 
\end{proposition}

This gives us a sufficient condition for 
the non-existence of a surjective homomorphism. 
By checking the divisibility of the Alexander polynomial, 
we prove the non-existence for many pairs of knots. 

Next, we apply the twisted Alexander polynomial to the remaining pairs 
which cannot be proved the non-existence by the Alexander polynomial. 
In this paper, 
we make use of $2$-dimensional unimodular representations of knot groups 
over finite prime fields. 
The twisted Alexander polynomial $\Delta_{K,\rho} (t)$, 
due to the Wada's definition, of $K$ 
for a representation $\rho : G(K) \to SL(2;\F_p)$ 
is defined to be a rational expression of one variable $t$ over $\F_p$, 
where $\F_p$ is the finite prime field of a characteristic $p$.  
See \cite{Wa1, KSW, KS1, KS2, KS3, KS4} 
for the precise definition and properties. 

Proposition \ref{alex-div} is extended to the following 
for the twisted Alexander polynomial. 

\begin{proposition}[Kitano-Suzuki-Wada \cite{KSW}]\label{ksw}
Assume that there exists a surjective homomorphism 
$\varphi : G(K_1) \to G(K_2)$. 
Then for any representation 
$\rho_2 : G(K_2) \to SL(2 ; \F_p)$ 
and $\rho_1 = \rho_2 \circ \varphi$, 
$\Delta_{K_1,\rho_1} (t)$ is divisible by $\Delta_{K_2,\rho_2} (t)$. 
More precisely, 
the denominator of $\Delta_{K_1,\rho_1} (t)$ is same as 
the denominator of $\Delta_{K_2,\rho_2} (t)$ and 
the numerator of $\Delta_{K_1,\rho_1} (t)$ can be divided by 
the numerator of $\Delta_{K_2,\rho_2} (t)$. 
\end{proposition}

Finally we can prove Proposition \ref{prop-nonex} 
by applying Proposition \ref{ksw}. 
The authors code two programs of Mathematica and Java independently 
and obtain Table \ref{non-ex}, 
which is used to prove Proposition \ref{prop-nonex}. 

This completes the proof of Theorem \ref{mainthm}. 

\section{Problems}\label{problems}

Here we mention some problems related to the main result. 

\leftline{(1) Simon's conjecture.}

The following question is arisen naturally from Theorem \ref{mainthm}. 

If there exists a surjective homomorphism
from $G(K_1)$ onto $G(K_2)$, 
then is the crossing number of $K_1$ 
greater than that of $K_2$?

Theorem \ref{mainthm} follows the affirmative answer 
in case the crossing number is smaller than or equal to $11$.
It still remains open for higher crossing cases. 

If the answer is affirmative in general, 
it turns out that Simon's conjecture \cite{K} holds, 
which is a part of problems on surjective homomorphisms between knot groups. 
In particular, the following is called Simon's conjecture: \\
For a given knot $K$, there exist only finitely many knot groups $G$ 
for which there is a surjective homomorphism 
$G(K) \rightarrow G$.

\leftline{(2) $2$-bridge knots.}

Recently Ohtsuki-Riley-Sakuma \cite{ORS} gave 
a systematic construction of surjective homomorphisms 
between $2$-bridge link groups. 
%
%
As shown in Table \ref{surj-31}-\ref{surj-63}, 
we constructed surjective homomorphisms explicitly. 
It is a problem 
whether they are same with Ohtsuki-Riley-Sakuma's construction. 

\begin{remark}
In \cite{BBRW}, 
it is announced that Simon's conjecture holds for any 2-bridge knot. 
\end{remark}

\leftline{(3) Geometric meanings. }

As we mention in the introduction, 
there are some reasons of the existence of surjective homomorphisms. 
In particular, it is important to determine which surjective homomorphism 
is induced by a degree one map. 
In other words, 
we should study 
characterization of surjective homomorphisms 
between knot groups induced by degree one maps.
 
\section{Tables}\label{table}

Now we describe how to read the tables below. 

We always take a Wirtinger presentation of $G(K)$ of a knot $K$: 
\[
G(K)= \langle x_1,\ldots,x_n |  r_1,\ldots, r_{n-1} \rangle. 
\]
For simplicity, 
we write a number representing a generator of $G(K)$. 
For example, we write 
$1,2,\dots,10,11$ 
for the generators 
$x_1 , x_2 , \ldots , x_{10} , x_{11}$ and 
$12\bar{1}\bar {10}$ means a relator 
$x_1 x_2 {x}_1^{-1} {x}_{10}^{-1}$. 
Here we fix the following presentations of the knot groups, 
which appear in the ranges of surjective homomorphisms. 
\begin{eqnarray*}
&&G(3_1) = 
\langle 
1,2,3 \, | \, 3 1 \bar{3} \bar{2}, \, 1 2 \bar{1} \bar{3}
\rangle ,\\
&&G(4_1) = 
\langle
1,2,3,4  \, | \, 
4 2 \bar{4} \bar{1},\, 1 2 \bar{1} \bar{3},\, 2 4 \bar{2} \bar{3}
\rangle ,\\
&&G(5_1) = 
\langle
1,2,3,4,5  \, | \, 
4 1 \bar{4} \bar{2},\, 5 2 \bar{5} \bar{3},\, 
1 3 \bar{1} \bar{4},\, 2 4 \bar{2} \bar{5} 
\rangle ,\\
&&G(5_2) = 
\langle
1,2,3,4,5
 \, | \,
4 1 \bar{4} \bar{2},\, 5 2 \bar{5} \bar{3},\, 
2 3 \bar{2} \bar{4},\, 1 4 \bar{1} \bar{5} 
\rangle ,\\
&&G(6_1) = 
\langle
1,2,3,4,5,6 \, | \,
4 2 \bar{4} \bar{1},\, 6 2 \bar{6} \bar{3},\, 5 3 \bar{5} \bar{4},\, 
1 5 \bar{1} \bar{4},\, 3 5 \bar{3} \bar{6}
\rangle ,\\
&&G(6_2) = 
\langle
1,2,3,4,5,6 \, | \,
4 2 \bar{4} \bar{1},\, 5 2 \bar{5} \bar{3},\, 6 3 \bar{6} \bar{4},\, 
1 5 \bar{1} \bar{4},\, 3 5 \bar{3} \bar{6}
\rangle ,\\
&&G(6_3) = 
\langle
1,2,3,4,5,6 \, | \,
3 2 \bar{3} \bar{1},\, 5 3 \bar{5} \bar{2},\, 6 3 \bar{6} \bar{4},\, 
2 5 \bar{2} \bar{4},\, 1 5 \bar{1} \bar{6}
\rangle .
\end{eqnarray*}
Under these notations, 
we give surjective homomorphisms to prove Theorem \ref{mainthm} 
in Table \ref{surj-31}-\ref{surj-63}. 
For any pair of knots with up to $10$ crossings in Theorem \ref{mainthm}, 
see the tables in \cite{KS1}.

In Table \ref{non-ex}, 
we list the pairs of knots $K,K'$ 
that there does not exist a surjective homomorphism from $G(K)$ onto $G(K')$. 
For example, 
there are $3_1$ and some knots with a prime integer in the first row of Table \ref{non-ex}. 
Here $11a_{6}(3)$ means 
that the non-existence of a surjective homomorphism from $G(11a_{6})$ onto $G(3_1)$ 
is proved by applying 
the twisted Alexander polynomial of $SL(2; \F_3)$-representations 
in Proposition \ref{ksw}. 

\begin{remark}
For the pairs of knots 
that do not appear in Table \ref{surj-31}-\ref{non-ex}, 
we can show that there exists no surjective homomorphism between their knot groups 
by using the classical Alexander polynomial. 
\end{remark}

\begin{longtable}{|c|l|}
\caption{Surjective homomorphisms to $3_1$}\\

\hline 
\label{surj-31}
$K$ & relators \\
\cline{2-2} 
& surjective homomorphism to $3_1$ \\
\hline 
\endfirsthead

\multicolumn{2}{c}%
{{\tablename\ \thetable{} Continued from previous page}} \\
\hline 
$K$ & relators \\
\cline{2-2} 
& surjective homomorphism to $3_1$ \\
\hline 
\endhead

\hline 
\multicolumn{2}{|r|}{{Continued on next page}} \\ 
\hline
\endfoot

\endlastfoot

\hline \tabtopsp{2mm}%
$11a_{43}$ & 
$10 2 \bar{10} \bar{1}, 8 3 \bar{8} \bar{2}, 6 4 \bar{6} \bar{3}, 
11 5 \bar{11} \bar{4}, 3 6 \bar{3} \bar{5}, 5 7 \bar{5} \bar{6}, 
2 8 \bar{2} \bar{7}, 7 9 \bar{7} \bar{8}, 1 10 \bar{1} \bar{9}, 
9 11 \bar{9} \bar{10}$ \\
\cline{2-2} \tabtopsp{2mm}%
& $1 \mapsto 1,2 \mapsto 3,3 \mapsto 2,
4 \mapsto 1,5 \mapsto 1,6 \mapsto 3,
7 \mapsto 2,8 \mapsto 1, 9 \mapsto 3, 
10 \mapsto 2, 11 \mapsto 1$ \\
\hline \tabtopsp{2mm}%
$11a_{44}$ & 
$10 2 \bar{10} \bar{1}, 8 3 \bar{8} \bar{2}, 6 3 \bar{6} \bar{4}, 
11 4 \bar{11} \bar{5}, 3 5 \bar{3} \bar{6}, 4 6 \bar{4} \bar{7}, 
2 8 \bar{2} \bar{7}, 7 9 \bar{7} \bar{8}, 1 10 \bar{1} \bar{9}, 
9 11 \bar{9} \bar{10}$ \\
\cline{2-2} \tabtopsp{2mm}%
& $1 \mapsto 1,2 \mapsto 1,3 \mapsto 3,
4 \mapsto 1,5 \mapsto 1,6 \mapsto 2,
7 \mapsto 3,8 \mapsto 2, 9 \mapsto 1, 
10 \mapsto 1, 11 \mapsto 1$ \\
\hline \tabtopsp{2mm}%
$11a_{46}$ & 
$10 2 \bar{10} \bar{1}, 8 3 \bar{8} \bar{2}, 6 3 \bar{6} \bar{4}, 
11 4 \bar{11} \bar{5}, 4 5 \bar{4} \bar{6}, 3 6 \bar{3} \bar{7}, 
2 8 \bar{2} \bar{7}, 7 9 \bar{7} \bar{8}, 1 10 \bar{1} \bar{9}, 
9 11 \bar{9} \bar{10}$ \\
\cline{2-2} \tabtopsp{2mm}%
& $1 \mapsto 1,2 \mapsto 2,3 \mapsto 1,
4 \mapsto 1,5 \mapsto 1,6 \mapsto 1,
7 \mapsto 1,8 \mapsto 3, 9 \mapsto 2, 
10 \mapsto \bar{1} 2 1, 11 \mapsto 1$ \\
\hline \tabtopsp{2mm}%
$11a_{47}$ & 
$11 2 \bar{11} \bar{1}, 7 2 \bar{7} \bar{3}, 1 4 \bar{1} \bar{3}, 
10 4 \bar{10} \bar{5}, 8 6 \bar{8} \bar{5}, 2 6 \bar{2} \bar{7}, 
9 8 \bar{9} \bar{7}, 6 9 \bar{6} \bar{8}, 4 9 \bar{4} \bar{10}, 
5 10 \bar{5} \bar{11}$ \\
\cline{2-2} \tabtopsp{2mm}%
& $1 \mapsto 1,2 \mapsto 1,3 \mapsto 1,
4 \mapsto 1,5 \mapsto 2,6 \mapsto 1,
7 \mapsto 1,8 \mapsto 3, 9 \mapsto 2, 
10 \mapsto 3, 11 \mapsto 1$ \\
\hline 
\end{longtable}

\begin{longtable}{|c|l|}
\caption{Surjective homomorphisms to $4_1$}\\

\hline 
\label{surj-41}
$K$ & relators \\
\cline{2-2} 
& surjective homomorphism to $4_1$ \\
\hline 
\endfirsthead

\multicolumn{2}{c}%
{{\tablename\ \thetable{} Continued from previous page}} \\
\hline 
$K$ & relators \\
\cline{2-2} 
& surjective homomorphism to $4_1$ \\
\hline 
\endhead

\hline 
\multicolumn{2}{|r|}{{Continued on next page}} \\ 
\hline
\endfoot

\endlastfoot

\hline \tabtopsp{2mm}%
$11a_{5}$ & 
$10 \, 1 \, \bar{10} \, \bar{2}, 6 \, 2 \, \bar{6} \, \bar{3}, 2 \, 4 \, \bar{2} \, \bar{3}, 
9 \, 4 \, \bar{9} \, \bar{5}, 3 \, 5 \, \bar{3} \, \bar{6}, 4 \, 7 \, \bar{4} \, \bar{6}, 
11 \, 8 \, \bar{11} \, \bar{7}, 5 \, 8 \, \bar{5} \, \bar{9}, 8 \, 10 \, \bar{8} \, \bar{9}, 
1 \, 10 \, \bar{1} \, \bar{11}$ \\
\cline{2-2} \tabtopsp{2mm}%
& $1 \mapsto 1,2 \mapsto 1 2 \bar{3},3 \mapsto 3 3 \bar{1} 4 \bar{1},
4 \mapsto 3,5 \mapsto 3,$ \\
\tabtopsp{2mm}%
& $6 \mapsto 1 2 2 \bar{3} \bar{3},
7 \mapsto 1 2 \bar{3},8 \mapsto 3, 9 \mapsto 3, 
10 \mapsto 3, 11 \mapsto 1 3 \bar{1}$ \\
\hline \tabtopsp{2mm}%
$11a_{6}$ & 
$8 \, 2 \, \bar{8} \, \bar{1}, 4 \, 3 \, \bar{4} \, \bar{2}, 11 \, 3 \, \bar{11} \, \bar{4}, 
6 \, 5 \, \bar{6} \, \bar{4}, 9 \, 5 \, \bar{9} \, \bar{6}, 10 \, 7 \, \bar{10} \, \bar{6}, 
1 \, 8 \, \bar{1} \, \bar{7}, 7 \, 9 \, \bar{7} \, \bar{8}, 5 \, 9 \, \bar{5} \, \bar{10}, 
2 \, 11 \, \bar{2} \, \bar{10}$ \\
\cline{2-2} \tabtopsp{2mm}%
& $1 \mapsto 1,2 \mapsto 1,3 \mapsto 4,
4 \mapsto 3,5 \mapsto 2 1 \bar{2},6 \mapsto 1 2 \bar{3},
7 \mapsto 1,8 \mapsto 1, 9 \mapsto 1, 
10 \mapsto 3, 11 \mapsto 2$ \\
\hline \tabtopsp{2mm}%
$11a_{51}$ & 
$8 \, 1 \, \bar{8} \, \bar{2}, 4 \, 3 \, \bar{4} \, \bar{2}, 11 \, 3 \, \bar{11} \, \bar{4}, 
7 \, 5 \, \bar{7} \, \bar{4}, 10 \, 5 \, \bar{10} \, \bar{6}, 5 \, 7 \, \bar{5} \, \bar{6}, 
9 \, 7 \, \bar{9} \, \bar{8}, 2 \, 8 \, \bar{2} \, \bar{9}, 6 \, 9 \, \bar{6} \, \bar{10}, 
1 \, 11 \, \bar{1} \, \bar{10}$ \\
\cline{2-2} \tabtopsp{2mm}%
& $1 \mapsto 1,2 \mapsto 1,3 \mapsto 4,
4 \mapsto 3,5 \mapsto 2,6 \mapsto 2 1 \bar{2},
7 \mapsto 1,8 \mapsto 1, 9 \mapsto 1, 
10 \mapsto 3, 11 \mapsto 2$ \\
\hline 
\end{longtable}

\begin{longtable}{|c|l|}
\caption{Surjective homomorphisms to $5_1$}\\

\hline 
\label{surj-51}
$K$ & relators \\
\cline{2-2} 
& surjective homomorphism to $5_1$ \\
\hline 
\endfirsthead

\multicolumn{2}{c}%
{{\tablename\ \thetable{} Continued from previous page}} \\
\hline 
$K$ & relators \\
\cline{2-2} 
& surjective homomorphism to $5_1$ \\
\hline 
\endhead

\hline 
\multicolumn{2}{|r|}{{Continued on next page}} \\ 
\hline
\endfoot

\endlastfoot

\hline \tabtopsp{2mm}%
$11n_{78}$ & 
$4 \, 2 \, \bar{4} \, \bar{1}, 9 \, 3 \, \bar{9} \, \bar{2}, 1 \, 4 \, \bar{1} \, \bar{3}, 
6 \, 4 \, \bar{6} \, \bar{5}, 11 \, 5 \, \bar{11} \, \bar{6}, 9 \, 7 \, \bar{9} \, \bar{6}, 
6 \, 8 \, \bar{6} \, \bar{7}, 2 \, 9 \, \bar{2} \, \bar{8}, 7 \, 10 \, \bar{7} \, \bar{9}, 
5 \, 10 \, \bar{5} \, \bar{11}$ \\
\cline{2-2} \tabtopsp{2mm}%
& $1 \mapsto 1,2 \mapsto 5,3 \mapsto 4,
4 \mapsto 3,5 \mapsto \bar{2} 4 2,6 \mapsto 4,
7 \mapsto \bar{2} 4 2,8 \mapsto 3, 9 \mapsto 2, 
10 \mapsto \bar{2} 1 2, 11 \mapsto 2$ \\
\hline \tabtopsp{2mm}%
$11n_{148}$ & 
$11 \, 2 \, \bar{11} \, \bar{1}, 9 \, 3 \, \bar{9} \, \bar{2}, 10 \, 3 \, \bar{10} \, \bar{4}, 
11 \, 5 \, \bar{11} \, \bar{4}, 3 \, 5 \, \bar{3} \, \bar{6}, 1 \, 7 \, \bar{1} \, \bar{6}, 
4 \, 8 \, \bar{4} \, \bar{7}, 5 \, 8 \, \bar{5} \, \bar{9}, 6 \, 10 \, \bar{6} \, \bar{9}, 
8 \, 10 \, \bar{8} \, \bar{11}$ \\
\cline{2-2} \tabtopsp{2mm}%
& $1 \mapsto 1,2 \mapsto \bar{2} 1 2,3 \mapsto 4,
4 \mapsto 5,5 \mapsto 4,6 \mapsto 4,
7 \mapsto 3,8 \mapsto 2, 9 \mapsto \bar{2} 3 2, 
10 \mapsto 2, 11 \mapsto 2$ \\
\hline 
\end{longtable}

\begin{longtable}{|c|l|}
\caption{Surjective homomorphisms to $5_2$}\\

\hline 
\label{surj-52}
$K$ & relators \\
\cline{2-2} 
& surjective homomorphism to $5_2$ \\
\hline 
\endfirsthead

\multicolumn{2}{c}%
{{\tablename\ \thetable{} Continued from previous page}} \\
\hline 
$K$ & relators \\
\cline{2-2} 
& surjective homomorphism to $5_2$ \\
\hline 
\endhead

\hline 
\multicolumn{2}{|r|}{{Continued on next page}} \\ 
\hline
\endfoot

\endlastfoot

\hline \tabtopsp{2mm}%
$11n_{71}$ & 
$4 \, 2 \, \bar{4} \, \bar{1}, 1 \, 3 \, \bar{1} \, \bar{2}, 7 \, 4 \, \bar{7} \, \bar{3}, 
2 \, 5 \, \bar{2} \, \bar{4}, 11 \, 5 \, \bar{11} \, \bar{6}, 9 \, 7 \, \bar{9} \, \bar{6}, 
3 \, 8 \, \bar{3} \, \bar{7}, 6 \, 9 \, \bar{6} \, \bar{8}, 8 \, 10 \, \bar{8} \, \bar{9}, 
1 \, 10 \, \bar{1} \, \bar{11}$ \\
\cline{2-2} \tabtopsp{2mm}%
& $1 \mapsto 1,2 \mapsto 5,3 \mapsto 4,
4 \mapsto 3,5 \mapsto 2,6 \mapsto 3,
7 \mapsto 2,8 \mapsto 1, 9 \mapsto 5, 
10 \mapsto 4, 11 \mapsto 5$ \\
\hline \tabtopsp{2mm}%
$11n_{185}$ & 
$7 \, 2 \, \bar{7} \, \bar{1}, 5 \, 3 \, \bar{5} \, \bar{2}, 1 \, 4 \, \bar{1} \, \bar{3}, 
11 \, 4 \, \bar{11} \, \bar{5}, 10 \, 6 \, \bar{10} \, \bar{5}, 1 \, 7 \, \bar{1} \, \bar{6}, 
2 \, 8 \, \bar{2} \, \bar{7}, 4 \, 9 \, \bar{4} \, \bar{8}, 3 \, 9 \, \bar{3} \, \bar{10}, 
2 \, 11 \, \bar{2} \, \bar{10}$ \\
\cline{2-2} \tabtopsp{2mm}%
& $1 \mapsto 1,2 \mapsto 5,3 \mapsto 5,
4 \mapsto 4,5 \mapsto 5,6 \mapsto 1 3 \bar{1},
7 \mapsto 3,8 \mapsto 2, 9 \mapsto 1, 
10 \mapsto 1 2 \bar{1}, 11 \mapsto 1$ \\
\hline 
\end{longtable}

\begin{longtable}{|c|l|}
\caption{Surjective homomorphisms to $6_1$}\\

\hline 
\label{surj-61}
$K$ & relators \\
\cline{2-2} 
& surjective homomorphism to $6_1$ \\
\hline 
\endfirsthead

\multicolumn{2}{c}%
{{\tablename\ \thetable{} Continued from previous page}} \\
\hline 
$K$ & relators \\
\cline{2-2} 
& surjective homomorphism to $6_1$ \\
\hline 
\endhead

\hline 
\multicolumn{2}{|r|}{{Continued on next page}} \\ 
\hline
\endfoot

\endlastfoot

\hline \tabtopsp{2mm}%
$11a_{352}$ & 
$5 \, 2 \, \bar{5} \, \bar{1}, 11 \, 3 \, \bar{11} \, \bar{2}, 10 \, 3 \, \bar{10} \, \bar{4}, 
9 \, 4 \, \bar{9} \, \bar{5}, 8 \, 6 \, \bar{8} \, \bar{5}, 3 \, 7 \, \bar{3} \, \bar{6}, 
2 \, 7 \, \bar{2} \, \bar{8}, 1 \, 9 \, \bar{1} \, \bar{8}, 4 \, 9 \, \bar{4} \, \bar{10}, 
7 \, 11 \, \bar{7} \, \bar{10}$ \\
\cline{2-2} \tabtopsp{2mm}%
& $1 \mapsto 1,2 \mapsto 2,3 \mapsto 2,
4 \mapsto 3,5 \mapsto 4,6 \mapsto 4,
7 \mapsto \bar{2} 1 4,8 \mapsto 4, 9 \mapsto 5, 
10 \mapsto 6, 11 \mapsto 2$ \\
\hline 
\end{longtable}

\begin{longtable}{|c|l|}
\caption{Surjective homomorphisms to $6_2$}\\

\hline 
\label{surj-62}
$K$ & relators \\
\cline{2-2} 
& surjective homomorphism to $6_2$ \\
\hline 
\endfirsthead

\multicolumn{2}{c}%
{{\tablename\ \thetable{} Continued from previous page}} \\
\hline 
$K$ & relators \\
\cline{2-2} 
& surjective homomorphism to $6_2$ \\
\hline 
\endhead

\hline 
\multicolumn{2}{|r|}{{Continued on next page}} \\ 
\hline
\endfoot

\endlastfoot

\hline \tabtopsp{2mm}%
$11a_{351}$ & 
$5 \, 2 \, \bar{5} \, \bar{1}, 11 \, 3 \, \bar{11} \, \bar{2}, 9 \, 3 \, \bar{9} \, \bar{4}, 
10 \, 4 \, \bar{10} \, \bar{5}, 8 \, 6 \, \bar{8} \, \bar{5}, 3 \, 7 \, \bar{3} \, \bar{6}, 
2 \, 7 \, \bar{2} \, \bar{8}, 1 \, 9 \, \bar{1} \, \bar{8}, 4 \, 9 \, \bar{4} \, \bar{10}, 
7 \, 11 \, \bar{7} \, \bar{10}$ \\
\cline{2-2} \tabtopsp{2mm}%
& $1 \mapsto 1,2 \mapsto 2,3 \mapsto 2,
4 \mapsto 3,5 \mapsto 4,6 \mapsto 4,
7 \mapsto \bar{2} 1 4,8 \mapsto 4, 9 \mapsto 5, 
10 \mapsto 6, 11 \mapsto 2$ \\
\hline 
\end{longtable}

\begin{longtable}{|c|l|}
\caption{Surjective homomorphisms to $6_3$}\\

\hline 
\label{surj-63}
$K$ & relators \\
\cline{2-2} 
& surjective homomorphism to $6_3$ \\
\hline 
\endfirsthead

\multicolumn{2}{c}%
{{\tablename\ \thetable{} Continued from previous page}} \\
\hline 
$K$ & relators \\
\cline{2-2} 
& surjective homomorphism to $6_3$ \\
\hline 
\endhead

\hline 
\multicolumn{2}{|r|}{{Continued on next page}} \\ 
\hline
\endfoot

\endlastfoot

\hline \tabtopsp{2mm}%
$11a_{47}$ & 
$4 \, 2 \, \bar{4} \, \bar{1}, 1 \, 3 \, \bar{1} \, \bar{2}, 8 \, 3 \, \bar{8} \, \bar{4}, 
2 \, 5 \, \bar{2} \, \bar{4}, 11 \, 5 \, \bar{11} \, \bar{6}, 9 \, 7 \, \bar{9} \, \bar{6}, 
3 \, 7 \, \bar{3} \, \bar{8}, 10 \, 9 \, \bar{10} \, \bar{8}, 7 \, 10 \, \bar{7} \, \bar{9}, 
5 \, 10 \, \bar{5} \, \bar{11}$ \\
\cline{2-2} \tabtopsp{2mm}%
& $1 \mapsto 1,2 \mapsto 6,3 \mapsto 5,
4 \mapsto 4,5 \mapsto 3,6 \mapsto 4,
7 \mapsto 3,8 \mapsto 2, 9 \mapsto 6, 
10 \mapsto 2 \bar{3} 4, 11 \mapsto 6$ \\
\hline \tabtopsp{2mm}%
$11a_{239}$ & 
$7 \, 2 \, \bar{7} \, \bar{1}, 1 \, 3 \, \bar{1} \, \bar{2}, 6 \, 4 \, \bar{6} \, \bar{3}, 
8 \, 4 \, \bar{8} \, \bar{5}, 9 \, 6 \, \bar{9} \, \bar{5}, 10 \, 6 \, \bar{10} \, \bar{7}, 
2 \, 8 \, \bar{2} \, \bar{7}, 11 \, 9 \, \bar{11} \, \bar{8}, 3 \, 9 \, \bar{3} \, \bar{10}, 
4 \, 11 \, \bar{4} \, \bar{10}$ \\
\cline{2-2} \tabtopsp{2mm}%
& $1 \mapsto 1,2 \mapsto 6,3 \mapsto 5,
4 \mapsto 5,5 \mapsto \bar{1} 3 5,6 \mapsto 5,
7 \mapsto 4,8 \mapsto 3, 9 \mapsto 3, 
10 \mapsto 2, 11 \mapsto 3$ \\
\hline 
\end{longtable}

\begin{longtable}{|c|l|}
\caption{Non-existence of surjective homomorphism} \\

\hline 
\label{non-ex}
$K'$ & $K\, (p)$ \\
\endfirsthead

\multicolumn{2}{c}%
{{\tablename\ \thetable{} Continued from previous page}} \\
\hline 
\endhead

\hline 
\multicolumn{2}{|r|}{{Continued on next page}} \\ 
\hline
\endfoot

\endlastfoot

\hline
$3_1$ & 
$
11a_{6}(3),
11a_{78}(3),
11a_{93}(3),
11a_{130}(3),
11a_{132}(5),
11a_{146}(3),
11a_{152}(3),
11a_{179}(5),
$ \\ & $
11a_{184}(5),
11a_{237}(3),
11a_{239}(3),
11a_{243}(3),
11a_{248}(3),
11a_{271}(3),
11a_{273}(7),
11a_{283}(3),
$ \\ & $
11a_{290}(3),
11a_{335}(5),
11a_{346}(5),
11a_{347}(3),
11a_{353}(5),
11a_{354}(3),
11a_{365}(5),
11n_{65}(3),
$ \\ & $
11n_{66}(3),
11n_{98}(3),
11n_{99}(3),
11n_{122}(7),
11n_{125}(3),
11n_{143}(3),
11n_{148}(3),
11n_{149}(3),
$ \\ & $
11n_{158}(3),
11n_{176}(3)
$ \\
\hline
$4_1$ & 
$
11a_{29}(3),
11a_{67}(3),
11a_{104}(13),
11a_{112}(29),
11a_{157}(53),
11a_{168}(17),
11a_{213}(3),
11a_{220}(3),
$ \\ & $
11a_{223}(3),
11a_{264}(7),
11a_{305}(7),
11a_{350}(3),
11a_{351}(3),
11a_{352}(3),
11n_{66}(3),
11n_{84}(3),
$ \\ & $
11n_{85}(3),
11n_{162}(7),
11n_{164}(7),
11n_{175}(7),
11n_{185}(3)
$ \\
\hline
\end{longtable}

\bibliographystyle{amsplain}

\end{document}